\documentclass[letterpaper, 10 pt, conference]{ieeeconf}  

\IEEEoverridecommandlockouts                              

\usepackage[pdftex]{graphicx}

\let\labelindent\relax

\usepackage{amsmath}
\usepackage{array}
\usepackage{amsfonts}
\usepackage{times}
\usepackage{algorithm}
\usepackage[noend]{algpseudocode}
\usepackage{xcolor} 
\usepackage[normalem]{ulem}
\usepackage{lipsum}

\usepackage{amsmath,amssymb,amsbsy}
\usepackage{amsthm} 
\usepackage{color}
\usepackage{comment}
\usepackage{graphicx}
\usepackage{enumerate}
\usepackage{epstopdf}
\usepackage{easy-todo}

\newcommand\norm[1]{\left\lVert#1\right\rVert}

\newcommand{\spr}[1] {{\color{red} #1} }
\newcommand{\carlo}[1] {{\color{blue} #1} } 

\newtheorem{theorem}{Theorem}[section]
\newtheorem{lemma}[theorem]{Lemma}

\begin{document}

\title{\LARGE \bf Optimal Control of Velocity and Nonlocal Interactions in the Mean-Field Kuramoto Model}

\author{Carlo~Sinigaglia, Francesco~Braghin, and Spring~Berman
%
\thanks{This work was supported by the Italian Ministry of Education, University and Research (MIUR).}
\thanks{C. Sinigaglia and F. Braghin are with the Department of Mechanical Engineering, Politecnico di Milano, Milano, Italy, 20156 {\tt\small \{carlo.sinigaglia, francesco.braghin\}@polimi.it}. S. Berman is with the School for Engineering of Matter, Transport and Energy, ASU, Tempe, AZ 85287 {\tt\small \{spring.berman\}@asu.edu}.
}}

\maketitle
\thispagestyle{empty}
\pagestyle{empty}

\begin{abstract}
In this paper, we investigate how the self-synchronization property of a swarm of Kuramoto oscillators can be controlled and exploited to achieve target densities and target phase coherence. In the limit of an infinite number of oscillators, the collective dynamics of the agents' density is described by a mean-field model in the form of a nonlocal PDE, where the nonlocality arises from the synchronization mechanism. In this mean-field setting, we introduce two space-time dependent control inputs to affect the density of the oscillators: an angular velocity field that corresponds to a state feedback law for individual agents,
and a control parameter that modulates the strength of agent interactions over space and time, i.e., a multiplicative control with respect to the integral nonlocal term.
We frame the density tracking problem as a PDE-constrained optimization problem. The controlled synchronization and phase-locking are measured with classical polar order metrics. After establishing the mass conservation property of the mean-field model and bounds on its nonlocal term,
a system of first-order necessary conditions for optimality is recovered using a Lagrangian method. 
The optimality system, comprising a nonlocal PDE for the state dynamics equation, the respective nonlocal adjoint dynamics, and the Euler equation, is
solved iteratively following a standard Optimize-then-Discretize approach and an efficient numerical solver based on spectral methods. We demonstrate our approach for each of the two control inputs in simulation.
\end{abstract} 


\IEEEpeerreviewmaketitle

\section{Introduction}

The emergence of collective behaviors from local interactions is ubiquitous in nature. Large groups of individuals can achieve complex macroscopic dynamical patterns by exploiting local interactions. The exact mechanisms by which flocks of birds and schools of fish respond to neighbors' behavior are still being investigated,
although various models have been proposed in mathematical biology to achieve such macroscopic behaviors in simulation that resemble their natural counterparts in a qualitative way; see, e.g., the Vicsek model \cite{vicsek}. Inspired by such biological systems, the robotics community has devised control strategies for large-scale collectives of robots 
with the goal of mimicking the resilience, efficiency, and adaptive behavior of biological swarms \cite{bio_robots}.  
Under the assumption of identical dynamics among individuals and a sufficiently large number of agents, it is possible to model the macroscopic behavior of a swarm of agents using a mean-field model  \cite{Elamvazhuthi2020}, which describes the density of agents at a specific location in the state space at a certain time. 
Mean-field models in the form of Louiville~\cite{Pogo}, Fokker-Planck~\cite{kart_2}, and McKean-Vlasov~\cite{Chen2021} equations are common in statistical physics and computational biology; they provide a powerful framework for studying the behavior 
of large-scale dynamical systems and have been recently analyzed from a control-theoretic perspective.

The steady-state collective behavior of such large-scale interacting systems is a prolific area of research in mathematical physics, along with the 
characterization of phase transitions in these systems with respect to a set of physical parameters. Very recently, the emergence of spatially inhomogeneous chimera states has been investigated for nonlocally interacting chiral particles \cite{Kruk2021}, while the effect of non-reciprocal interactions between different species of individuals has been studied in \cite{Fruchart2021}.

In this paper, we investigate how interactions, synchronization patterns, and phase transitions can be exploited from a control perspective. We consider a large-scale swarm of identical Kuramoto oscillators, in which two possible control mechanisms affect
the dynamics of each individual oscillator:  
an angular velocity control input
and a control input that locally modulates the strength of oscillator interactions. 
We formulate an Optimal Control Problem (OCP) subject to the mean-field nonlocal Partial Differential Equation (PDE) model that describes the dynamics of the oscillator density, and we derive a 
set of first-order necessary conditions for optimality for this system. We then show in simulation how the control inputs 
can guide the density of oscillators to target distributions and accelerate their synchronization. 

\section{State Dynamics}
The Kuramoto oscillator network \cite{kuramoto2003chemical} is a classical model that describes 
synchronization in a variety of biological and engineering systems. It consists of an array of coupled oscillators that each affect the others' phases via pairwise interactions.
We consider a swarm of $N$ oscillators in which the time-dependent phase $\theta_i(t)$ of each oscillator $i$ is governed by the following controlled dynamical model:   
\begin{equation*}
    \dot{\theta}_i = \omega_i ~+~ u_1 ~+~ \frac{u_2}{N} \sum_{j=1}^{N} \sin(\theta_j - \theta_i -\alpha) \, + ~\eta(t), ~~~ i=1,...,N,
\end{equation*}
where $\omega_i$ is the oscillator's 
natural frequency, $\alpha$ is a constant phase shift that gives 
rise to traveling wave solutions  \cite{Kruk2021}, 
$\eta(t)$ is a white noise stochastic process whose intensity is associated with diffusion coefficient $D>0$, 
$u_1$ is an angular velocity control input, 
and $u_2$ is a control gain 
that modulates the strength of agent interactions. 
The usual Kuramoto phase model \cite{kuramoto2003chemical} is 
recovered by setting $u_1 = 0$ and $u_2=K$, where $K$ is a possibly non-constant interaction strength. Depending on the sign and magnitude of $K$ and its dependence on $\theta$, 
a variety of oscillator phase distributions, and hence degrees of synchronization, can appear \cite{Abrams2004,Hong2011}.

Kuramoto oscillators with noise are known to exhibit synchronization once the quantity $\frac{K}{D}$ exceeds a threshold,
assuming that $K>0$. Furthermore, traveling wave solutions and spatially inhomogeneous phases known as chimera states are known to arise for particular combinations of the interaction kernel 
and the phase shift $\alpha$ 
\cite{Hong2016,Hong2011,Abrams2004}. 
The complex polar order parameter defined below is used 
to measure phase coherence and synchronization intensity: 
\begin{equation}
\label{disc_polar}
    R(t)e^{i \psi(t)} = \frac{1}{N} \sum_{j=1}^{N} e^{i \theta_j(t)}
\end{equation}
The amplitude $R$ takes values in the interval $[0,1]$, where $0$ corresponds to disordered motion (incoherence) and $1$ to global synchronization (phase locking). The variable $\psi$ denotes the mean-field phase.

For the sake of simplicity, we assume that all oscillators have the same 
natural frequency, $\omega_i = \omega$. Without loss of generality, we can select $\omega=0$, which is equivalent to expressing the oscillator dynamics in a frame rotating at angular velocity $\omega$. 

In the limit of an infinite number of oscillators, i.e., $N \to \infty$, the oscillator density 
$q(\theta,t)$ satisfies the following mean-field nonlocal PDE,
\begin{equation}
\label{mf_state}
    q_t - D q_{\theta \theta} + \, \partial_{\theta} \left( u_2 \, w[q] q + u_1\, q \right) = 0, 
\end{equation}
where 
\begin{equation}
\label{wDef}
w[q] = \int_{S^{1'}} \sin(\theta'-\theta-\alpha) \, q(\theta',t) d\theta'
\end{equation} 
is a 
nonlocal integral term describing the oscillator interactions. Note that we use the notation $S^{1'}$ for the unit circle to highlight the role of the integration variable $\theta'$.

Eq. (\ref{mf_state}) is a nonlocal parabolic PDE 
that 
describes the evolution in space-time of the probability density of a single oscillator having phase $\theta$ at time $t$, or equivalently, the normalized concentration of a population of oscillators at this phase and time. 
The natural functional space for $q$ is the space of absolutely continuous probability measures, that is, $q \in V$ where $V=\{ f \in L^2(0,T,H^1(S^1)) : \int_{S^1} f(\theta,t) d\theta = 1 \, \textrm{and} \,  f \geq 0 \,\, a.e. \,\, t \in (0,T)\}$. In the following, we will refer to 
the solution of Eq.  (\ref{mf_state}) with $u_1 \equiv 0$ and $u_2 \equiv K$ as the {\it uncontrolled mean-field dynamics}. In this case, we will set $K=1$.

Equation (\ref{mf_state}) is defined on the unit circle $S^1$ and is thus naturally equipped with periodic boundary conditions. The mean-field equivalent \cite{Kruk2021} of the complex polar order parameter is:
\begin{equation}
\label{cont_polar}
    R(t)e^{i \psi(t)} = \int_{S^1} e^{i \theta} q(\theta,t) \, d\theta 
\end{equation}
This continuous version of the polar order parameter has the same range for $R$ as its discrete counterpart, since the oscillator density $q$ is normalized to have unitary total mass.
 
In the following, we will make frequent use of Green's theorem for functions defined on $S^1$. Given functions $f, g: S^1 \rightarrow \mathbb{R}$, we have that: 
\begin{equation}
\label{s_1_green}
    \int_{S^1} f_{\theta} g \, d\theta ~=~  -\int_{S^1} f g_{\theta} \, d\theta.
\end{equation}
We will also need the following identity, 
\begin{equation}
\label{s_1_conv}
    \int_{S^1} w[f] g \, d\theta ~=~  \int_{S^1} w^*[g] f \,  d\theta,
\end{equation}
where the functional $w^{*}$ is defined as:
\begin{equation*}
w^*[g] = \int_{S^{1'}} \sin(\theta-\theta'-\alpha) \, g(\theta') d\theta'.
\end{equation*}
It is straightforward to prove Eq. (\ref{s_1_conv}) 
by exchanging the variables of integration.

It is also easy to show that Eq. (\ref{mf_state}) conserves the total mass, $m(t) = \int_{S^1} \, q \, d\theta$, for any choice of control functions $u_1$, $u_2$ and parameters $D$, $\alpha$. First, note that Eq. (\ref{mf_state}) can be equivalently written as follows, where $-Dq_{\theta \theta}$ is written as $\partial_{\theta}(-D q_{\theta})$: 
\begin{equation}
\label{div_state}
    q_t + \, \partial_{\theta} \left( u_2 \, w[q] q + u_1\, q -D q_{\theta} \right) = 0.
\end{equation} 
Then, we have that
\begin{equation*}
    \dot{m} = \int_{S^1} q_t \, d\theta 
     = - \int_{S^1} \, \partial_{\theta} \left( u_2 \, w[q] q + u_1\, q -D q_{\theta}\right) d\theta = 0,
\end{equation*}
where we have substituted the state dynamics (\ref{div_state}) and 
applied Eq. (\ref{s_1_green}) with $g\equiv1$. Furthermore, the velocity $w[q]$ of the transport field $w[q]\,q$, which is induced by the \emph{nonlocal} wave, is bounded by 1, as we demonstrate in the next lemma. 
\begin{lemma}
\label{w_bound}
For all probability densities $q(t) \in H^1(S^1)$, 
    $\norm{w[q(t)]}_{L^{\infty}(S^1)} \leq 1$ for a.e. $t \in (0,T)$ 
\begin{proof}
Using Eq. \eqref{wDef} and the fact that $q$ is a probability density, we obtain the following inequality:
\begin{equation*}
\begin{aligned}
    \left|w[q(t)]\right| &~=~ \left|\int_{S^{1'}} \sin(\theta'-\theta-\alpha) \, q(\theta',t) d\theta' \right| \\ 
    &~\leq~ \sup_{\theta'}\left(|\sin(\theta'-\theta-\alpha)|\right) \int_{S^{1'}} 
    |q| d\theta' ~\leq~ 1.
\end{aligned}
\end{equation*}
Consequently, we have that
\begin{equation*}
     \norm{w[q(t)]}_{L^{\infty}(S^1)} ~=~ \sup_{q(t) \in H^1(S^1)}|w[q(t]| ~\leq~ 1.
\end{equation*}
\end{proof}
\end{lemma}
In the following lemma, we use Lemma \ref{w_bound} to obtain an upper bound on
the norm of the mean-field state $q$.


\begin{lemma}
Defining the set of control inputs as $\mathcal{U}=L^{2}(0,T,L^{\infty}(S^1))$, that is, $u_1,u_2 \in \mathcal{U}$, it can be shown that
\begin{equation}
\label{normqBound}
\norm{q}_V^2 ~\leq~C_{\mathcal{U}}(u_1,u_2)\norm{q_0}_{L^2(S^1)}^2,
\end{equation}
where $q_0 \in L^2(S^1)$ is the initial probability density and 
$C_{\mathcal{U}}(u_1,u_2) \equiv \norm{u_1}_{\mathcal{U}} 
+ \norm{u_2}_{\mathcal{U}}$. 

\begin{proof}
We just need to show that G\aa rding's inequality holds for this problem, and then we can proceed as in Theorem 1 in \cite{sinigaglia2021density}. The weak form of Eq. (\ref{mf_state}) can be written 
as:
\begin{equation*}
    \int_{S^1} q_t(t)\phi\,d\theta +  a(q(t),\phi) = 0 \hspace{5mm} \forall \phi \in H^1(S^1)
\end{equation*}
for a.e. $t\in(0,T)$, where the form $a(q,\phi)$, defined as 
\begin{equation*}
    a(q,\phi) = \int_{S^1} ( D q_{\theta} \phi_{\theta} - (u_2 w[q] +u_1)q \, \phi_{\theta} ) \, d\theta,
\end{equation*}
is obtained by multiplying Eq. (\ref{mf_state}) by a test function $\phi \in H^1(S^1)$, integrating the equation over 
the domain $S^1$, and then applying Eq. (\ref{s_1_green}). We need to show that for almost every $t \in (0,T)$, there exists $\lambda(t) > 0$ such that:  
\begin{equation*}
    a(q(t),q(t)) + \lambda(t)\norm{q(t)}_{L^2(S^1)}^2 ~\geq~ \alpha_0(t) \norm{q(t)}_{H^1(S^1)}^2
\end{equation*}
for some $\alpha_0(t) \geq 0$; i.e., that the form $a$ is weakly coercive. First, by applying Lemma \ref{w_bound}, the Cauchy-Schwarz inequality, and Cauchy's inequality with $\epsilon$
(see e.g. \cite{evans}, Appendix B.2), we can show that for every $\epsilon \geq 0$, 
\begin{equation*}
\begin{aligned}
    & \left| \, \int_{S^1} (u_2 w[q] +u_1)q \, q_{\theta} d\theta \, \right| \\ 
    &~~\leq \norm{(u_2 w[q] +u_1)}_{L^{\infty}(S^1)} \norm{q}_{L^2(S^1)} \norm{q_{\theta}}_{L^2(S^1)} \\
    & ~~\leq \Big( \norm{u_{2}}_{L^{\infty}(S^1)} + \norm{u_{1}}_{L^{\infty}(S^1)} \Big)\norm{q}_{L^2(S^1)} \norm{q_{\theta}}_{L^2(S^1)} \\
    & ~~ = C
    \norm{q}_{L^2(S^1)} \norm{q_{\theta}}_{L^2(S^1)} \\
    & ~~\leq   \frac{C^2} 
    {4\epsilon}\norm{q}_{L^2(S^1)}^2 + \epsilon  \norm{q_{\theta}}^2_{L^2(S^1)},
\end{aligned} 
\end{equation*}
where $C \equiv C_{L^{\infty}(S^1)}(u_1,u_2) \equiv \norm{u_1}_{L^{\infty}(S^1)} + \norm{u_2}_{L^{\infty}(S^1)}$.
Therefore, for a.e. $t\in(0,T)$, we can now write: 
\begin{equation*}
\begin{aligned}
    & a(q,q) + \lambda \norm{q}_{L^2(S^1)}^2 \\
    &\geq   D\norm{q_{\theta}}^2_{L^2(S^1)} +  \lambda \norm{q}_{L^2(S^1)}^2 \\
    &~~~~~ - \left|\int_{S^1} (u_2 w[q] +u_1)q \, q_{\theta} d\theta \right| \\
    & \geq (D-\epsilon)\norm{q_{\theta}}^2_{L^2(S^1)} +\Big(\lambda-\frac{C^2}{4\epsilon}\Big)\norm{q}_{L^2(S^1)}^2 \\
    & \geq \frac{D}{2}\norm{q_{\theta}}^2_{L^2(S^1)} +\Big(\lambda-\frac{C^2}{2D}\Big)\norm{q}_{L^2(S^1)}^2 \\
    & \geq \frac{D}{2}\norm{q_{\theta}}^2_{L^2(S^1)} +\frac{C^2}{2D}\norm{q}_{L^2(S^1)}^2 \\  
    & \geq \alpha_0\norm{q}^2_{H^1(S^1)},  \\
\end{aligned}
\end{equation*}
where we have chosen $\epsilon = \frac{D}{2}$, 
$\lambda= \frac{C^2}{D}$, and $\alpha_0 = \min\{\epsilon,\frac{\lambda}{2}\}$. Note that the time dependence of variables in the inequalities above is omitted to improve readability. 

We can thus 
conclude that the form $a$ is weakly coercive and apply the same reasoning as in  \cite{sinigaglia2021density}, Theorem 1, to establish 
a bound on the norm of the mean-field state 
for every $u_1 , u_2 \in \mathcal{U}$. 
\end{proof}
\end{lemma}
\section{The Optimal Control Problem}\label{Section:ControlDesign}

In this section, we formulate an OCP subject to the mean-field dynamics 
in Eq. (\ref{mf_state}). A natural objective is to compute the minimum-energy control inputs $u_1$, $u_2$ that drive the oscillator density to a target density
$z(\theta,t)$ at a 
final time $t=T$. We define the cost functional to be minimized as 
\begin{equation*}
J(q,u_1,u_2) = J_q(q) + J_u(u_1,u_2),
\end{equation*}
where
\begin{equation*}
\begin{aligned}
J_q(q) = ~& \frac{\alpha_r}{2} \int_{0}^T \int_{S^1} (q(\theta,t)-z(\theta,t))^2 \, d \theta \, dt \\
&+ \frac{\alpha_t}{2} \int_{S^1} (q(\theta,T)-z(\theta,T))^2 d\theta, 
\end{aligned}
\end{equation*}
\begin{equation*}
J_u(u_1,u_2) = \frac{1}{2} \int_{0}^T\int_{S^1} (\beta_1 u_1(\theta,t)^2 + \beta_2 u_2(\theta,t)^2) \\    d \theta dt,
\end{equation*}
and $\alpha_r$, $\alpha_t$, $\beta_1$, and $\beta_2$ are nonnegative weighting constants.

The OCP can be written as:
\begin{equation}
\label{ocp_formulation}
\begin{aligned}
&\min_{u_1,u_2,q} \quad  J(q,u_1,u_2) =  J_q(q) +  J_u(u_1,u_2)                
   \\
&\text { s.t. } \begin{array}{ll}
\displaystyle     q_t - D q_{\theta \theta} + \, \partial_{\theta}\left( u_2 \, w[q] q + u_1\, q \right) = 0,   
\phantom{space} \\
\phantom{space} \\
q(\theta,0) = q_0(\theta) .\\
\end{array}
\end{aligned}
\end{equation}
We consider 
the control inputs $u_1$ and $u_2$ to be 
acting simultaneously in the OCP, and we derive a set of first-order necessary optimality conditions accordingly using a Lagrangian method. 
In order to do so, we define the Lagrangian functional \cite{Trltzsch2010OptimalCO} as:
\begin{equation*}
\begin{aligned}
    &\mathcal{L}(q,u,p) = J(q,u_1,u_2) \\
    & ~ - \int_{0}^{T}\int_{S^1} \left( q_t - D p_{\theta \theta}+ \, \partial_{\theta} \left( u_2 \, w[q] q + u_1\, q \right) \right) p \, \, d\theta dt,
\end{aligned}
\end{equation*}
where $p: S^1 \times (0,T) \mapsto \mathbb{R}$ is the adjoint field still to be determined. The adjoint dynamics are recovered by taking the Gateaux derivative with respect to the state $q$ and setting it to zero for any state variation. We focus on deriving the first variation of the third term of the PDE constraint in (\ref{ocp_formulation}), since the treatment of the time derivative and diffusion terms 
is standard in the literature on OCPs 
for parabolic equations (see e.g. \cite{Trltzsch2010OptimalCO}). 
First, it is useful to note that: 
\begin{equation}
\label{conv_term}
\begin{aligned}
    & \int_{0}^{T}\int_{S^1} \, \partial_{\theta}  \Big(u_2 \, w[q] (q)\Big)p \, \, d\theta dt \\
    & ~~~ = ~ - \int_{0}^{T}\int_{S^1} \,u_2 \, w[q] q \, p_{\theta} \,  d\theta dt.
\end{aligned}
\end{equation}
Taking the first variation of the term of interest, we have:
\begin{equation*}
    \begin{aligned}
          &\frac{\partial}{\partial \epsilon}\Big\rvert_{0}\int_{S^1} -\,u_2 \, w[q+\epsilon \psi] (q+\epsilon \psi) \, p_{\theta} \, \, d\theta dt \\
          &~~~~ = \int_{S^1}-u_2 \Big(w[q]\psi + w[\psi]q\Big)p_{\theta} \, \, d\theta dt \\
          &~~~~ = -\int_{S^1}\Big(u_2 w[q] p_{\theta} +w^{*}[u_2 p_{\theta}\,q]\Big) \psi \, \, d\theta dt,
    \end{aligned}
\end{equation*}
where we have used Eq. (\ref{s_1_conv}) to take out $\psi$ from the nonlocal operator $w$.
Regarding the component of the term that contains 
$u_1$, we use Eq. (\ref{s_1_green}) to write:
\begin{equation*}
    \int_{S^1} \partial_{\theta} \Big( \, u_1 \, q \Big) p \, d\theta ~=~ - \int_{S^1} \, u_1 \, q \, p_{\theta} d \theta 
\end{equation*}
so that the state variation of this term reads:
\begin{equation}
\label{flux}
\begin{aligned}
    \frac{\partial}{\partial \epsilon}\Big\rvert_{0}  \int_{S^1} \, u_1 \, \Big(q+\epsilon \psi \Big) p_{\theta} d \theta ~~=~ \int_{S^1} \, u_1 \, p_{\theta} \psi d \theta. 
\end{aligned}
\end{equation}

The adjoint dynamics can be then written as: 
\begin{equation}
\label{adj_dyn}
\displaystyle     
-p_t - D p_{\theta \theta} - \Big(u_2 \,  w[q] + u_1 \Big) p_{\theta}  
 - w^{*}[u_2\,p_{\theta}\,q] = \alpha_r \left(q - z\right),  
\end{equation}
with the final time condition specified as:
\begin{equation*}
    p(\theta,T) = \alpha_t \left( q(\theta,T)-z(\theta,T) \right).
\end{equation*}

We can now derive the Euler equation for the reduced gradient with respect to the controls $u_1$ and $u_2$. Using Eq. (\ref{flux}), the first variation of the Lagrangian with respect to 
the velocity control input $u_1$ is straightforward to compute: 
\begin{equation}
\label{grad1}
    \nabla J_{u_1}(\theta,t) = \beta_1 u_1(\theta,t) + q(\theta,t) \, p_{\theta}(\theta,t).
\end{equation}
This control mechanism for mean-field PDEs,  
i.e., optimal control of a velocity field, has been considered in various applications (see e.g. \cite{sinigaglia2021density,kart_2}), and controllability results \cite{karthik_bil} have been obtained for the case where agent interactions are not present (i.e., $u_2 \equiv 0$).
In order to recover the reduced gradient with respect to $u_2$, we make use of Eq. (\ref{conv_term}), which is linear in $u_2$. We can then find that 
the reduced gradient has the form: 
\begin{equation}
\label{grad2}
    \nabla J_{u_2}(\theta,t) = \beta_2 u_2(\theta,t) + w[q](\theta,t)\,q(\theta,t)\,p_{\theta}(\theta,t).
\end{equation}
Note that the integral term $w[q]$ makes this gradient equation nonlocal, in contrast to gradient equation (\ref{grad1}), which means that 
at each point in space-time, the 
equation 
depends on the entire 
state solution $q$. 
 
\section{Simulation Results and Discussion}\label{Section:Simulation}

In this section, we solve the OCP numerically using an iterative gradient descent method and analyze the results.  
At each iteration, the current control input values $u_1$, $u_2$ are used to solve the state and adjoint equations, Eqs. (\ref{mf_state}) and (\ref{adj_dyn}). The reduced gradient is then computed using Eqs. (\ref{grad1}) and (\ref{grad2}), and the control inputs are updated such that they produce a 
sufficient decrease 
in the cost functional. The numerical algorithm is implemented in Python and the operators are discretized using the open-source software Dedalus \cite{Burns2020}, which makes efficient use of novel spectral algorithms for the solution of PDEs. Although we have derived optimality conditions for the general case where both $u_1$ and $u_2$ are applied, we optimize $u_1$ and $u_2$ separately in the simulations
in order to investigate the ability of each input to steer the state dynamics towards the target density.



In the {\it$u_1$-controlled} case, we 
set $u_2 \equiv K = 1$ and 
design the velocity control input $u_1(\theta,t)$ to steer
an initial density of oscillators towards a target density $z(\theta)$ at final time $T$. 
A similar control problem is solved in \cite{kart_2}, which does not include agent interactions, and in \cite{Chen2021}, which includes interactions in the form of spatially attractive/repulsive potentials.
Our aim is to investigate the effect of the control input $u_1$ on the synchronization mechanism, as compared to the uncontrolled dynamics ($u_1 = 0, u_2 = 1$).

In the {\it $u_2$-controlled case}, we set $u_1 = 0$ and design the interaction strength control input $u_2(\theta,t)$, which locally regulates the intensity of agent interactions, to solve this density steering problem.
At each iteration of the optimization algorithm, 
$u_2(\theta,t)$ is updated using Eq. (\ref{grad2}); 
the adjoint equation \eqref{adj_dyn} is the same for both the $u_1$-controlled and $u_2$-controlled cases. 
We note that the influence of the integral term $w[q]$ may limit the control authority with respect to 
$u_2$; these restrictions on controllability must be further investigated. 

In all simulations, we set $T=10$ s and choose 
$D=0.25$ and $\alpha = 0$, for which a synchronized phase is the solution of the uncontrolled system at steady-state \cite{Kruk2021} in the form of a Gaussian density that can be characterized semi-analytically. The objective of the OCP 
is to speed up the convergence of the state dynamics to the synchronized phase while steering 
the mean phase to the mean of the target density. Toward this end, we define the target density as a relatively low-variance Gaussian function with a mean of $\frac{3}{2}\pi$ rad.   
 
\begin{figure}[h]
 \centering
 \includegraphics[width =  \linewidth]{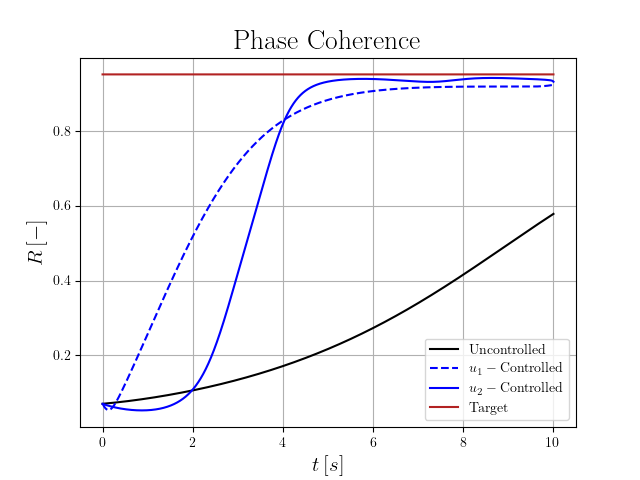}
  \caption{Phase coherence parameter 
  $R(t)$ over time. Both the velocity and interaction strength control inputs
  speed up the system self-synchronization.
  }  
\label{R}
 \end{figure}

\begin{figure}[h]
 \centering
 \includegraphics[width =  \linewidth]{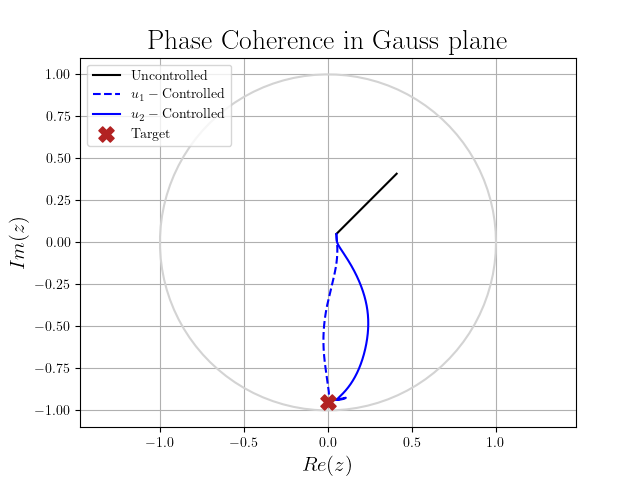}
   \caption{Image of the continuous synchronization metric $R(t)e^{i\psi(t)}:[0,T]\mapsto \mathbb{C}$. Trajectories of both controlled systems converge to the target coherence and mean phase.}
\label{psi}
 \end{figure}

\begin{figure}[h]
 \centering
 \includegraphics[width =  \linewidth]{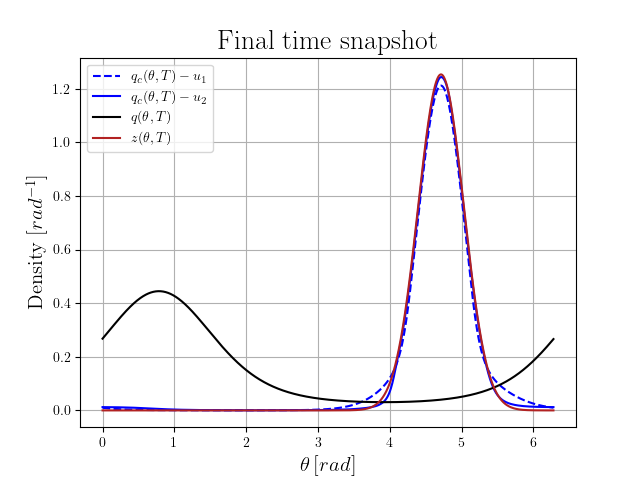}
 \caption{Snapshot of the state $q(\theta,T)$ 
 at final time $T=10$ s.  
 The uncontrolled system {\it (black)} 
has not reached a coherent phase, while the controlled systems {\it (blue)} are both close to the target density {\it (red)}.
 }
\label{T}
 \end{figure}

Fig. \ref{R} plots the time evolution of the phase coherence parameter, the amplitude of the complex polar order parameter in Eq. (\ref{cont_polar}), for both controlled systems and the uncontrolled system.  The trajectories of the phase coherence parameter and mean phase in the complex plane for all three systems are plotted
in Fig. \ref{psi}, along with the target coherence and mean phase. These figures show that both control inputs $u_1$ and $u_2$ are able to steer the system towards the target coherence and mean phase, with $u_1$ producing synchronization more quickly than $u_2$. In contrast, the uncontrolled system achieves synchronization much more slowly, and it convergences to a mean phase that differs from the target phase. Figure \ref{T} plots the $u_1$-controlled and $u_2$-controlled densities and the uncontrolled density at the final time $T = 10$ s, along with the target density. The controlled densities have both closely approached the target density, while the uncontrolled system has not yet synchronized. 

\begin{figure}[h]
 \centering
 \includegraphics[width =  \linewidth]{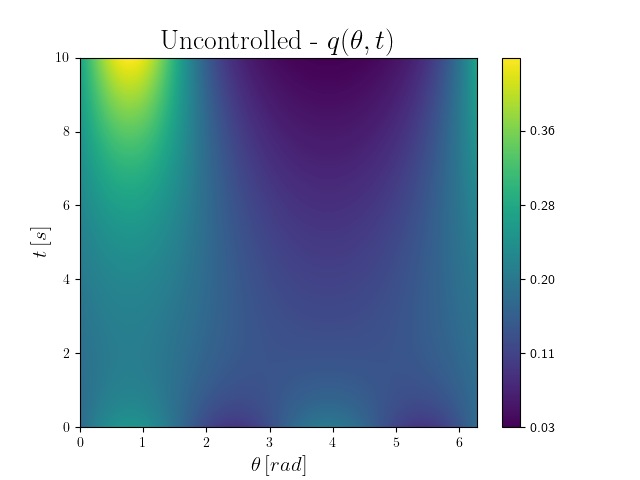}
 \caption{Uncontrolled state dynamics with diffusion parameter $D=0.25$ and 
 interaction strength $K=1$. Synchronization emerges close to the final time, with mean phase dependent on the initial conditions. }
\label{tc_1_u_unc}
 \end{figure}
\begin{figure}[h]
 \centering
 \includegraphics[width =  \linewidth]{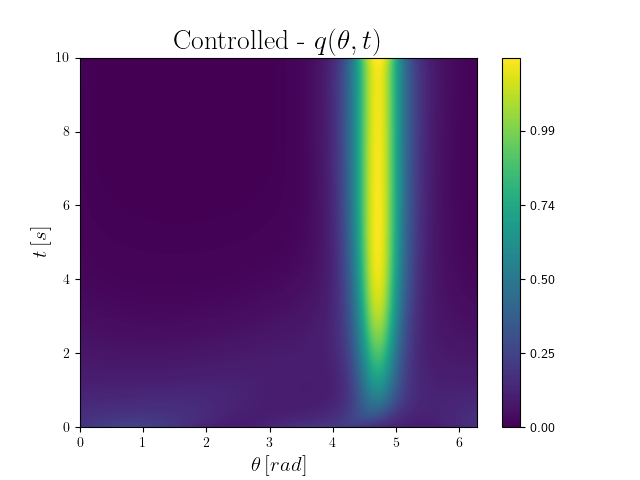}
 \caption{Controlled state dynamics under optimal velocity control input 
 $u_1$. The density is coherent and synchronized after a brief transient period. 
 }
\label{tc_1_u_opt}
 \end{figure}
\begin{figure}[h]
 \centering
 \includegraphics[width =  \linewidth]{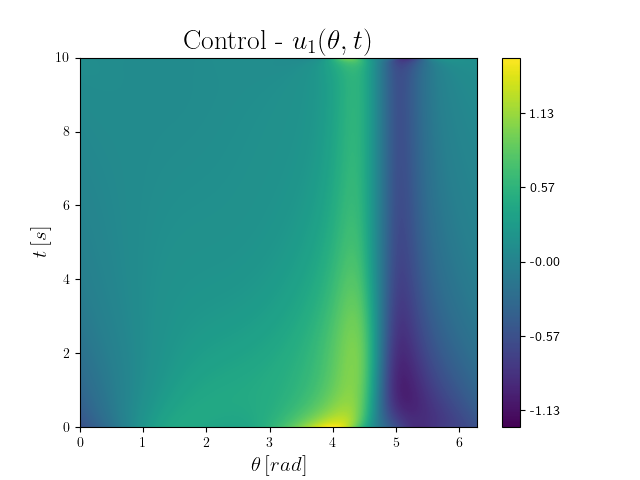}
 \caption{Space-time evolution of the optimal velocity 
control input $u_1(\theta,t)$. 
 Once the state $q(\theta,t)$  is steered close to 
 the target density, the control  magnitude is highest in the vicinity of the mean phase of the target density, 
 inducing both positive and negative local transport fields.}
\label{tc_1_C}
 \end{figure}
\begin{figure}[h]
 \centering
 \includegraphics[width =  \linewidth]{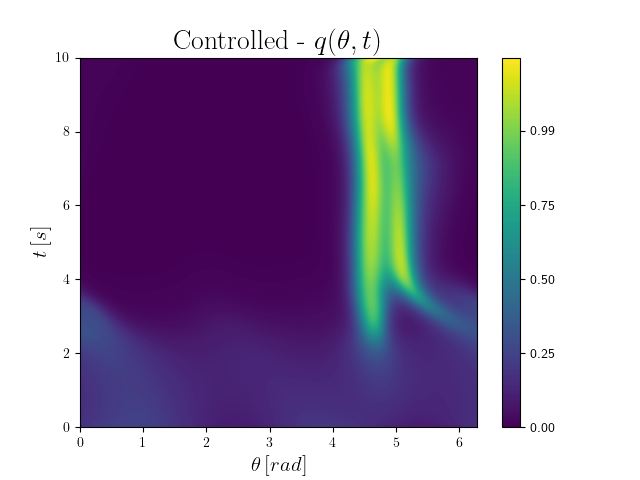}
 \caption{Controlled state dynamics under optimal interaction strength control input 
 $u_2$. The transient period before synchronization accelerates 
 is longer than in the $u_1$-controlled case, 
 and the steady-state density is less coherent.}  
\label{tc_2_q}
 \end{figure}
\begin{figure}[h]
 \centering
 \includegraphics[width =  \linewidth]{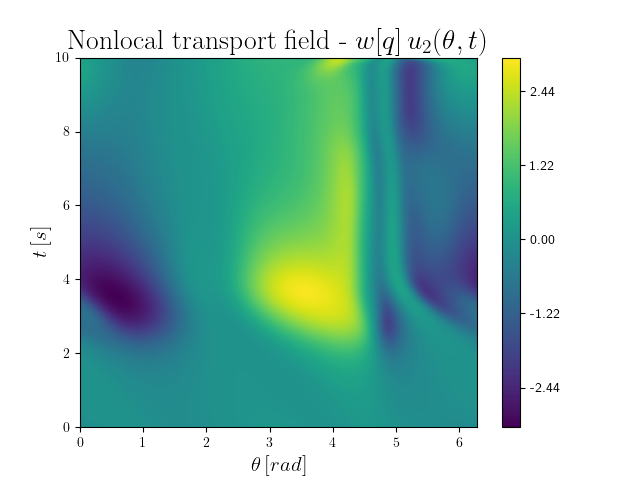}
 \caption{Space-time evolution of the controlled nonlocal transport field $w[q]u_2$ induced by the optimal interaction strength control input $u_2(\theta,t)$. 
 A negative control action is associated with 
 oscillators that tend to disalign. The control intensity, given by $u_2$ (not shown), is approximately seven times higher than $u_1$ in the $u_1$-controlled case due to the influence of the nonlocal term $w[q]$, which directly multiplies $u_2$. 
 }
\label{tc_2_c} 
 \end{figure}

The state dynamics of the uncontrolled, $u_1$-controlled, and $u_2$-controlled systems over the entire simulation are shown in Figs. \ref{tc_1_u_unc}, \ref{tc_1_u_opt}, and \ref{tc_2_q}, respectively. The $u_2$-controlled system exhibits more complex dynamics than the $u_1$-controlled system before converging to the target density; this is due to the previously mentioned limitations on control authority with respect to $u_2$. The space-time evolution of the corresponding control input $u_1$ and the controlled nonlocal transport field $w[q]u_2$ are plotted in Figs. \ref{tc_1_C} and \ref{tc_2_c}, respectively. The quantity $w[q]u_2$ is plotted rather than $u_2$ in order to directly compare the controlled transport fields produced by $u_1$ and $u_2$.
 

 
 
It is worthwhile to comment on the signs of the control inputs. 
The positive (negative) sign of the velocity control input $u_1$ 
generates a transport field in the direction of increasing (decreasing) $\theta$, while the role of the sign for the interaction strength control input $u_2$ is more subtle. 
Positive values of $u_2$ drive the corresponding oscillator density toward alignment, and thus synchronization, while negative values produce oscillator disalignment
\cite{Hong2011}. Thus, in order to reach a 
target density, both locally repelling and locally aligning interaction strength inputs  are needed. 

In summary, both the velocity and interaction strength control inputs are able to successfully drive the oscillator density to a target distribution while increasing the self-synchronization speed.   
Interaction strength control is able to achieve similar performance as velocity control, with slightly slower convergence to synchronization, despite its controllability limitations. 

\section{Conclusion} 

In this paper, we have presented 
a novel Optimal Control Problem (OCP) for a large-scale system of identical Kuramoto oscillators in the mean-field limit. The mean-field state dynamics in the OCP 
are nonlocal due to an integral term describing oscillator interactions; 
this term is shown to be a bounded functional, and some preliminary results on 
the well-posedness of the state dynamics are established. Both control mechanisms considered, the angular velocity and interaction strength of the oscillators,
are able to accurately track a target density, as well as increase the synchronization speed over that of the uncontrolled dynamics. 
Furthermore, the interaction strength control input achieves the target density by 
exploiting the nonlocal term. These 
results set the grounds for establishing controllability and existence results for OCPs in which 
the state dynamics are governed by self-synchronizing nonlocal PDEs, which in turn are envisioned to find application in the control of collective systems such as robotic swarms through these 
synchronization mechanisms. 

\bibliographystyle{IEEEtran}
\bibliography{References.bib}
\end{document}